\documentclass[runningheads,a4paper]{llncs}

\usepackage[T1]{fontenc}
\usepackage{amssymb}
\usepackage{dsfont}
\usepackage{amsmath}
\usepackage{graphicx}

\allowdisplaybreaks

\newtheorem{dfn}{Definition}

\newtheorem{thm}{Theorem}

\newtheorem{rmk}{Remark}
\newtheorem{exm}{Example}

\title{Prediction of quantiles by statistical learning and
application to GDP forecasting}

\titlerunning{Prediction of quantiles and GDP forecasting}

\author{Pierre Alquier\inst{1,3} and Xiaoyin Li \inst{2}}
\institute{LPMA (Universit\'e Paris 7)
 \\ 175, rue du Chevaleret \\ 75013 Paris FRANCE \\
\email{alquier@math.jussieu.fr}
\\
\email{http://alquier.ensae.net/}
\\
\and Laboratoire de Math\'ematiques (Universit\'e de Cergy-Pontoise) \\
  UCP site Saint-Martin, 2 boulevard Adolphe Chauvin \\ 95000 Cergy-Pontoise, FRANCE \\
\email{xiaoyin.li@u-cergy.fr}
\\
\and CREST (ENSAE)}
\date{\today}

\authorrunning{Pierre Alquier and Xiaoyin Li}

\toctitle{Prediction of quantiles by statistical learning and
application to GDP forecasting}
\tocauthor{Pierre Alquier and Xiaoyin Li}

\begin{document}

\maketitle

\begin{abstract}
In this paper, we tackle the problem of prediction and confidence intervals
for time series using a statistical learning approach and quantile loss
functions. In a first time, we show that the Gibbs estimator is able to predict
as well as the
best predictor in a given family for a wide set of loss functions. In particular,
using the quantile loss function of \cite{quantileReg}, this
allows to build confidence intervals. We apply these results to the
problem of prediction and confidence regions for the French Gross Domestic Product
(GDP) growth, with promising results.
\end{abstract}

\section{Introduction}

Motivated by economics problems, the prediction of time series is one of the
most emblematic problems of statistics. Various methodologies are used that come
from such various fields as parametric statistics, statistical learning, computer
science or game theory.

In the parametric approach, one assumes that the time series is generated from
a parametric model, e.g. ARMA or ARIMA, see \cite{hamilton,davis}.
It is then possible to estimate the parameters of the model and to build confidence
intervals on the prevision. However, such an assumption is unrealistic in most applications.

In the statistical learning point of view, one usually tries to avoid such restrictive
parametric assumptions - see, e.g., \cite{lugosi,stoltz} for the online approach
dedicated to the prediction of individual sequences, and \cite{modha,meir,alqwin}
for the batch approach. However, in this setting, a few attention has been paid to
the construction of confidence intervals or to any quantification of the precision
of the prediction. This is a major drawback in many applications.
Notice however that \cite{quantileBiau} proposed to minimize the cumulative risk
corresponding to the quantile loss function defined by \cite{quantileReg}. This
led to asymptotically correct confidence intervals.

In this paper, we propose to adapt this approach to the batch setting and provide
nonasymptotic results. We also apply these results to build quarterly prediction and
confidence regions for the French Gross Domestic Product (GDP) growth. Our approach
is the following. We assume that we are given a set of basic predictors - this is
a usual approach in statistical learning, the predictors are sometimes referred
as ``experts'', e.g. \cite{lugosi}. Following \cite{alqwin}, we describe a procedure
of aggregation, usually referred as Exponentially Weigthed Agregate (EWA),
\cite{DT1,gerchi}, or Gibbs estimator, \cite{Catoni2004,Catoni2007}.
It is interesting to note that this procedure is also related to aggregations
procedure in online learning as the weighted majority algorithm of \cite{LiWa},
see also \cite{VOVK}. We give a
PAC-Bayesian inequality that ensures optimality properties for this procedure.
In a few words, this inequality claims that our predictor performs as well as the
best basic predictor up to a remainder of the order $\mathcal{K}/\sqrt{n}$ where
$n$ is the number of observations and $\mathcal{K}$ measures the complexity of
the set of basic predictors. This result is very general, two conditions will
be required: the time series must be weakly dependent in a sense that we will
make more precise in Section~\ref{section_theorem}, and
loss function must be Lipschitz.
This includes, in particular, the quantile loss functions. This allows us
to apply this result to our problem of economic forecasting.

The paper is organized as follows:
Section~\ref{section_context} provides notations used in the whole paper.
We give a definition of the Gibbs estimator in Section~\ref{section_description}.
The PAC-Bayesian inequality is given in
Section~\ref{section_theorem}, and the application to GDP forecasting
in Section~\ref{section_application}.
Finally, the proof of all the theorems are given in
Section~\ref{sectionproofs}.

\section{The context}
\label{section_context}

Let us assume that we observe
$X_{1},\ldots,X_{n}$ from a $\mathbb{R}^{p}$-valued stationary time series
$X=\left(X_{t}\right)_{t\in\mathbb{Z}}$ defined on
$(\Omega,\mathcal{A},\mathbb{P})$.
From now, $\|.\|$ will denote the Euclidian norm on $\mathbb{R}^{p}$.
Fix an integer $k$ and let us assume that we are given a 
family of predictors $\left\{
f_{\theta}:(\mathbb{R}^{p})^{k}\rightarrow\mathbb{R}^p,\theta\in\Theta\right\}$:
 for any $\theta$ and any $t$, $f_{\theta}$ applied to the last past values
$\left(X_{t-1},\ldots,X_{t-k}\right)$ is a possible
prediction of $X_{t}$. For the sake of simplicity, let us put for any
$t\in\mathbb{Z}$ and any $\theta\in\Theta$,
$$ \hat{X}_{t}^{\theta} = f_{\theta}(X_{t-1},\ldots, X_{t-k}).$$
We also assume that $\theta\mapsto f_{\theta}$ is linear.
Note that we may want to include parametric models as well as non-parametric
prediction. In order to deal with
various family of predictors, we propose a model-selection type approach:
$ \Theta = \cup_{j=1}^{m} \Theta_{j} $.
\begin{exm}
We deal with only one model, $m=1$ and $\Theta=\Theta_{1}$.
We put
$\theta=(\theta_0,\theta_1,\ldots,\theta_k)
\in\Theta = \mathbb{R}^{k+1}$ and define the linear autoregressive predictors
$$ f_{\theta}(X_{t-1},\ldots, X_{t-k}) = \theta_0 + \sum_{j=1}^{k} \theta_{j}
X_{t-j} . $$
\end{exm}
\begin{exm}
We may generalize the previous example to non-parametric auto-regression, for
example using a dictionary of functions
$ (\mathbb{R}^{p})^{k} \rightarrow\mathbb{R}^{p} $, say
$(\varphi_{i})_{i=0}^{\infty}$. Then we can fix $m=n$,
and take
$\theta=(\theta_1,\ldots,\theta_\ell)\in\Theta_{j}=\mathbb{R}^{j}$ and
$$ f_{\theta}(X_{t-1},\ldots, X_{t-k}) = \sum_{i=1}^{j}
\theta_i \varphi_{i}(X_{t-1},\ldots,X_{t-k}) . $$
\end{exm}
We now define a quantitative criterion to evaluate the quality
of the predictions. Let $\ell$ be a loss function, we will assume that $\ell$
satisfies
the following assumption.

\noindent {\bf Assumption LipLoss}: $\ell$ is given by: $\ell(x,x')=g(x-x')$ for
some convex function $g$ satisfying $g\geq 0$, $g(0)=0$ and $g$ is $K$-Lipschitz.
\begin{dfn} We put, for any $\theta\in\Theta$, $
R\left(\theta\right)=\mathbb{E}\left[\ell\left(\hat{X}_{t}^{\theta},X_{t}
\right)\right]$.
\end{dfn}
Note that because of the stationarity, $R(\theta)$ does not depend on $t$.
\begin{exm}
A first example is $\ell(x,x')=\|x-x'\|$. In this case, the Lipschitz constant
$K$ is $1$.
This example was studied in detail in \cite{alqwin}.
In \cite{modha,meir}, the loss function is the quadratic loss 
$\ell(x,x')=\|x-x'\|^2$. Note that it also satisfies our Lipschitz condition,
but only if we assume that the time series is bounded.
\end{exm}
\begin{exm}
\label{exm_quantile}
When the time-series is real-valued, we can use a quantile loss function.
The class of quantile loss functions is defined as
$$
\ell_{\tau}(x,y)=
\begin{cases}
\tau\left( x-y\right), & \text{if } x-y>0\\
-\left(1-\tau\right)\left( x-y\right), & \text{otherwise}
\end{cases}
$$
where $\tau\in\left(0,1\right)$.
It is motivated by the
following remark: if $U$ is a real-valued random variable, then any value
$t^*$ satisfying $ \mathbb{P}(U\leq t^*) = \tau$ is a minimizer of
of $t\mapsto \mathbb{E}(\ell_{\tau}(X-t))$; such a value is called
quantile of order $\tau$ of $U$.
This loss function was introduced by \cite{quantileReg}, see \cite{quantileBook}
for a survey. Recently, \cite{cherno} used it in the context of high-dimensional
regression, and \cite{quantileBiau} in learning problems.
\end{exm}

\section{Gibbs estimator}
\label{section_description}

In order to introduce the Gibbs estimator, we first define the empirical risk.

\begin{dfn}
For any $\theta\in\Theta$,
$ r_{n}(\theta) = \frac{1}{n-k} \sum_{i=k+1}^{n} \ell \left(
\hat{X}_{i}^{\theta},X_{i}\right) .$
\end{dfn}

Let $\mathcal{T}$ be a $\sigma$-algebra on $\Theta$ and $\mathcal{T}_{\ell}$ be
its restriction to $\Theta_{\ell}$. Let $\mathcal{M}_{+}^{1}(\Theta)$
denote the set of all probability
measures on $(\Theta,\mathcal{T})$, and $\pi\in\mathcal{M}_{+}^{1}(\Theta)$.
This probability measure is usually called the {\it prior}. It will
be used to control the complexity of the set of predictors $\Theta$.

\begin{rmk}
In the case where $\Theta=\cup_{j}\Theta_{j}$ and the $\Theta_{j}$ are disjoint,
we can write
$ \pi({\rm d}\theta ) = \sum_{j = 1}^{m} \mu_{j} \pi_{j}({\rm
d}\theta) $
where $\mu_{j}:=\pi(\Theta_{j})$ and $\pi_{j}({\rm d}\theta):=\pi({\rm
d}\theta) \mathbf{1}_{\Theta_{j}}(\theta)
/\mu_{j}$. Here $\pi_{j}$ can be interpreted as a prior probability
measure inside the model $\Theta_{j}$
and that the $\mu_{j}$ as a prior probability measure between the models.
\end{rmk}

\begin{dfn}[Gibbs estimator]
\label{def_est}
 We put, for any $\lambda>0$,
$$ \hat{\theta}_{\lambda} = \int_{\Theta} \theta \hat{\rho}_{\lambda}({\rm
d}\theta),\text{ where } \hat{\rho}_{\lambda}({\rm d}\theta) = \frac{e^{-\lambda
r_{n}(\theta)}\pi({\rm d}\theta) }
                                    { \int e^{-\lambda r_{n}(\theta')}\pi({\rm
d}\theta') }. $$
\end{dfn}
 The choice of the parameter $\lambda$ is discussed in the next sections.

\section{Theoretical results}
\label{section_theorem}

In this section, we provide a PAC-Bayesian oracle inequality for the Gibbs
estimator. PAC-Bayesian bounds were introduced in \cite{STW97,McA2},
see \cite{Catoni2004,Catoni2007,AlquierPAC,AudibertHDR,CatoniAudibert}
for more recent advances. The idea is that the risk of the
Gibbs estimator will be close to $\inf_{\theta}R(\theta)$ up to a small remainder.
More precisely, we upper-bound it by
$$ \inf_{\rho\in\mathcal{M}_{+}^{1}(\Theta)}
\left\{ \int R(\theta)\rho({\rm d}\theta) + \text{ remainder}(\rho,\pi)
\right\}. $$
To establish such a result,
we need some hypothesis. The first hypothesis concerns the type of dependence of
the process, it uses the
$\theta_{\infty,n}(1)$-coefficients of \cite{Dedecker2007a}. Such a condition
is also
used in \cite{alqwin}, and is more general than the mixing conditions
used in \cite{meir,modha}.

\noindent {\bf Assumption WeakDep}: we assume that the distribution $\mathbb{P}$
is such that almost surely,
$\|X_0\|_{\infty}\leq \mathcal{B} <\infty$, and that there is
a constant $\mathcal{C}$ with $\theta_{\infty,k}(1)\leq \mathcal{C} < \infty$
for any $k$. We remind that
for any $\sigma$-algebra $\mathfrak{S}\subset\mathcal{A}$, for any
$q\in\mathbb{N}$,
for any $(\mathbb{R}^{p})^{q}$-valued random variable $Z$ defined on
$(\Omega,\mathcal{A},\mathbb{P})$, we define
$$ \theta_{\infty}(\mathfrak{S},Z) = \sup_{f\in\Lambda_{1}^{q}} \Bigl\|
\mathbb{E}\left[f(Z)|\mathfrak{S}\right]
- \mathbb{E}\left[f(Z)\right] \Bigr\|_{\infty} $$
where
$$\Lambda_{1}^{q} = \left\{f:(\mathbb{R}^{p})^{q} \rightarrow \mathbb{R}, \quad
                   \frac{|f(z_1,\ldots,z_q) -
f(z'_1,\ldots,z'_q)|}{\sum_{j=1}^{q} \|z_j - z'_j\| } \leq 1 \right\} ,$$
and
$$\theta_{\infty,k}(1) := \sup \left\{ \theta_{\infty}(\sigma(X_t,t\leq
p),(X_{j_1},\ldots,X_{j_\ell})) ,
      \quad p < j_1 < \ldots < j_\ell, 1 \leq \ell \leq k \right\}.  $$
Intuitively, these coefficients provide a quantification of the dependence
between the past and the future of the time series. The sequence $\theta_{\infty,k}(1)$
is growing with $k$, but, when $X_k$ behaves almost like a random variable independent
from $X_0$, $X_{-1}$, $...$, the sequence is bounded.
Examples of processes satisfying {\bf WeakDep} are provided
in \cite{alqwin,Dedecker2007a}. It includes processes
of the form $ X_t = H(\xi_t,\xi_{t-1},\xi_{t-2},\dots) $
where the $\xi_{t}$ are iid and bounded and $H$ satisfies a Lipschitz
condition, in particular, ARMA processes with bounded innovations.
It also includes uniform $\varphi$-mixing processes (see \cite{Doukhan1994})
and some dynamical systems.

\noindent {\bf Assumption Lip}: we assume that there is a constant $L>0$ such that
for any $\theta\in\Theta$, there
are coefficients
$a_{j}\left(\theta\right)$ for $1\leq j \leq k$ satisfying, for any $x_1$, ...,
$x_k$ and $y_1$, ..., $y_k$,
$$
\left\|f_{\theta}\left(x_1,\ldots,x_k\right)-f_{\theta}\left(y_1,\ldots,
y_k\right)\right\|
     \leq \sum_{j=1}^{k} a_j\left(\theta\right)\left\Vert x_j-y_j\right\Vert,
\text{ with } \sum_{j=1}^{k}a_{j}\left(\theta\right) \leq L.
$$

\begin{thm}[PAC-Bayesian Oracle Inequality]
\label{main_result}
Let us assume that assumptions {\bf LipLoss}, {\bf WeakDep} and {\bf Lip} are
satisfied. Then, for any $\lambda>0$, for any $\varepsilon>0$, with probability
at least $1-\varepsilon$,
$$ R\left(\hat{\theta}_{\lambda}\right)  \leq
\inf_{\rho\in\mathcal{M}_{+}^{1}(\Theta)} \left[ \int R(\theta) \rho({\rm d}\theta)
 + \frac{2\lambda \kappa^{2}}{n\left(1-\frac{k}{n}\right)^2} +
\frac{2\mathcal{K}(\rho,\pi) +  2 \log\left(\frac{2}{\varepsilon}\right)
}{\lambda} \right]
$$
where $\kappa = \kappa(K,L,\mathcal{B},\mathcal{C}) := K(1+L) (\mathcal{B} +
\mathcal{C})/\sqrt{2}$ and $\mathcal{K}$ is the
Kullback divergence, given by $ \mathcal{K}(\rho,\pi)=
 \int \log \left[\frac{{\rm d}\rho}{{\rm d}\pi}(\theta)\right] \rho({\rm d}\theta)$
if $\rho \ll \pi$ and $+\infty$ otherwise.
\end{thm}

The choice of $\lambda$ in practice is a hard problem. In
\cite{Catoni2003,Catoni2007} a general method is proposed
to optimize the bound with respect to $\lambda$. However, while
adapted in the iid case, this method is more difficult to use in
the context of time series as it would require the knowledge of $\kappa$.
However,
some empirical calibration seems to give
good results, as shown in Section \ref{section_application}.

At the price of a more technical analysis,
this result can be extended to the case where the $X_t$ are not assumed
to be bounded: the results in~\cite{alqwin} require subGaussian tails
for $X_t$, but suffer a $\log(n)$ loss in the learning rate.

\section{Application to French GDP and quantile prediction}

\label{section_application}

\subsection{Uncertainty in GDP forecasting}

Every quarter $t$, economic forecasters at INSEE\footnote{{\it Institut
National de la Statistique et des Etudes Economiques}, the French national bureau of
statistics, http://www.insee.fr/} are asked a prediction for
the quarterly growth rate of the French Gross Domestic Product (GDP). Since it involve
a lot of information, the ``true value'' of the growth rate $\log({\rm GDP}_t
/ {\rm GDP}_{t-1})$ is only known after two years,
but {\it flash estimates} of the growth rate, say $\Delta {\rm GDP}_{t}$, are
published 45 days after the end of the current quarter $t$.
One of the most relevant economic pieces of information available at
time $t$ to the forecaster, apart from past GDP
observations, are {\it business surveys}. Indeed, they are a rich source of information,
for at least two reasons. First, they are rapidly available, on a monthly basis.
Moreover, they provide information coming directly from the true economic decision
makers.

A business survey is traditionally a fixed questionnaire of ten questions sent monthly
to a panel of companies. This process is described in \cite{Devilliers}. INSEE publishes
a composite indicator called the {\it French business climate indicator}: it summarises
information of the whole survey. This indicator is defined in \cite{climate}, see also
\cite{DuboisMichaux}. All these values are available from the INSEE website.
Note that a quite similar approach is used in other countries, see also
\cite{predBiau} for a prediction of the European Union GDP based on EUROSTATS data.

It is however well known among economic forecasters that interval confidence or
density forecasts are to be given with the prediction, in order to provide an
idea of the uncertainty of the prediction. The ASA and the NBER started using
density forecasts in 1968, see \cite{Diebold,Tay} for historical surveys on density
forecasting. The Central Bank of England and INSEE,
among others, provide their prediction with a ``fan chart'', \cite{Britton}.
However, it is interesting to note that the methdology used is often very crude,
see the criticism in \cite{CornecCIRET,Dowd}. For example, until 2012, the fan
chart provided by the INSEE led to the construction of confidence intervals with
constant length. But there is an empirical evidence that it is more difficult to
forecast economic quantities during crisis (e.g. the subprime crisis in 2008).
The Central Bank of
England fan chart is not reproducible as it includes subjective information.
Recently, \cite{CornecCIRET} proposed a clever density forecasting method based on
quantile regressions that gives satisfying results in practice.
However, this method did not receive any theoretical support up to our knowledge.
 
Here, we use the Gibbs estimator to
build a forecasting of $\Delta {\rm GDP}_{t}$, using the quantile loss function.
This allows to return a prediction: the forecasted
median, for $\tau=0.5$, that is theoretically supported. This also allows to provide
confidence intervals corresponding to various quantiles.

\subsection{Application of Theorem \ref{main_result}}

At each quarter $t$, the objective is to predict the flash estimate of GDP growth,
$\Delta {\rm GDP}_{t}$. As described previouly, the available information is $\Delta {\rm GDP}_{t'}$
for $t'<t$ and $I_{t'}$ for $t'<t$, where for notational convenience, $I_{t-1}$ is the climate
indicator available to the INSEE at time $t$ (it is the mean of the climate indicator at month
3 of quarter $t-1$ and at month 1 and 2 of quarter $t$). The observation period is
1988-Q1 (1st quarter of 1988) to 2011-Q3.

We define $X_t=(\Delta {\rm GDP}_{t},I_t)'\in\mathbb{R}^{2}$.
As we are not interested by the prevision of $I_t$ but only by the prediction
of the GDP growth, the loss function will only take into account $\Delta {\rm
GDP}_{t}$. We use the quantile
loss function of Example \ref{exm_quantile}:
\begin{multline*}
\ell_{\tau}((\Delta {\rm GDP}_{t},I_t),(\Delta' {\rm GDP}_{t},I'_t))
\\
=
\begin{cases}
\tau\left( \Delta {\rm GDP}_{t} - \Delta' {\rm GDP}_{t} \right), & \text{if }
\Delta {\rm GDP}_{t} - \Delta' {\rm GDP}_{t} >0\\
-\left(1-\tau\right)\left( \Delta {\rm GDP}_{t} - \Delta' {\rm GDP}_{t} 
\right), & \text{otherwise}.
\end{cases}
\end{multline*}

To keep in mind that the risk depends on $\tau$, we add a subscript $\tau$ in the notation
$ R^{\tau}(\theta) := \mathbb{E} \left[ \ell_{\tau} \left( \Delta {\rm
GDP}_{t}, f_{\theta}(X_{t-1},X_{t-2})\right) \right]$.
We also let $r_{n}^{\tau}$ denote the associated empirical risk.
Following \cite{CornecCIRET,stage} we consider predictors of the form:
\begin{equation}
\label{modele}
f_{\theta}(X_{t-1},X_{t-2})=\theta_0+\theta_1 \Delta {\rm GDP}_{t-1} +
\theta_{2} I_{t-1} + \theta_{3}
                      (I_{t-1}-I_{t-2})|I_{t-1} - I_{t-2}|
\end{equation}
where $\theta=(\theta_0,\theta_1,\theta_2,\theta_3)\in \Theta(B)$.
For any $B>0$ we define
$$\Theta(B)=
\biggl\{\theta=(\theta_0,\theta_1,\theta_2,\theta_3)\in\mathbb{R}^{4},
 \|\theta\|_1=\sum_{i=0}^3 |\theta_i| \leq B \biggr\} .$$
These predictors of Equation~\ref{modele} correspond to the model used in
\cite{CornecCIRET} for forecasting,
one of the conclusions of \cite{CornecCIRET,stage} is that these family of predictors
allow to obtain a forecasting as precise as the INSEE one.

For technical reason that will become clear in the proofs, if one wants to achieve
a prediction performance
comparable to the best $\theta\in\Theta(B)$, it is more convenient to define the
prior $\pi$ as the uniform
probability distribution on some slightly larger set, e.g. $\Theta(B+1)$. We will
let $\Pi_B$ denote this
distribution. We let
$\hat{\rho}^{\tau}_{B,\lambda}$ and $\hat{\theta}_{B,\lambda}^{\tau}$
denote repectively the associated agregation distribution and the associated
estimator, defined in Definition~\ref{def_est}.

In this framework, Assumption {\bf Lip} is satisfied with $L = B+1$,
and
the loss function is $K$-Lipschitz with $K=1$ so Assumption {\bf LipLoss} is also
satisfied.

\begin{thm}
\label{thm_appli}
Let us fix $\tau\in(0,1)$.
Let us assume that Assumption {\bf WeakDep} is satisfied, and that $n\geq
\max\left(10,\kappa^2 / (3\mathcal{B}^{2})\right)$.
Let us fix $\lambda = \sqrt{3n}/\kappa$. Then, with probability at least
$1-\varepsilon$ we have
\begin{equation*}
 R^{\tau} (\hat{\theta}_{B,\lambda}^{\tau}) 
\leq \inf_{\theta\in\Theta(B)}\left\{ R^{\tau}(\theta) +
\frac{2\sqrt{3} \kappa}{\sqrt{n}} \left[
2.25 + \log\left(\frac{(B+1)\mathcal{B} \sqrt{n}}{\kappa}\right)
+ \frac{\log\left(\frac{1}{\varepsilon}\right)}{3}
   \right]\right\}.
\end{equation*}
\end {thm}
The choice of $\lambda$ proposed in the theorem may be a problem as in practice
we will not know $\kappa$.
Note that from the proof, it is obvious that in any case, for $n$ large enough,
when $\lambda=\sqrt{n}$ we still have a bound
\begin{equation*}
 R^{\tau} (\hat{\theta}_{B,\lambda}^{\tau}) 
\leq \inf_{\theta\in\Theta(B)}\left\{ R^{\tau}(\theta) +
\frac{ C(B,\mathcal{B},\kappa,\varepsilon)}{\sqrt{n}} \right\}.
\end{equation*}
However, in practice, we will work in an online setting: at each date $t$
we compute the Gibbs estimator based on the observations from $1$ to $t$ and use
it to predict the GDP and its quantiles at time $t+1$. Let
$\hat{\theta}^{\tau}_{B,\lambda}[t]$ denote this estimator. We propose the following
empirical approach: we define a set of values $\Lambda=\{2^{k},k\in\mathbb{N}\}\cap\{1,...,n\}$.
At each step $t$, we compute $\hat{\theta}^{\tau}_{B,\lambda}[t]$ for each $\lambda\in\Lambda$
and use for prediction $ \hat{\theta}^{\tau}_{B,\lambda(t)}[t]$ where
$\lambda(t)$ is defined by
$$ \lambda(t) = \arg\min_{\lambda\in\Lambda} \sum_{j=3}^{t-1} \ell_{\tau}(\Delta GDP_{j},
f_{\hat{\theta}^{\tau}_{B,\lambda}[j]}(X_{j-1},X_{j-2})), $$
namely, the value that is currently the best for online prediction. This choice
leads to good numerical results.

In practice, the choice of $B$ has less importance. As soon as $B$ is large enough,
simulation shows that the estimator
does not really depend on $B$, only the theoretical bound does.
As a consequence we take $B=100$ in our experiments.

\subsection{Implementation}

We use the importance sampling method to compute
$\hat{\theta}_{B,\lambda}^{\tau}[t]$, see, e.g., \cite{Christian}.
We draw an iid sample $T_1$, ..., $T_N$ of vectors in $\mathbb{R}^{4}$, from the
distribution
$\mathcal{N}(\hat{\theta}^{\tau},vI)$ where $v>0$ and $\hat{\theta}^{\tau}$ is simply the
$\tau$-quantile regression estimator of $\theta$ in~\eqref{modele}, as computed by the
``quantile regression package'' of the R software \cite{R}. Let $g(\cdot)$ denote the density
of this
distribution. By the law of large numbers we have
$$ 
         \sum_{i=1}^{N}   \frac{ T_i
\exp\left[-\lambda r_{t}(T_i) \right]\mathbf{1}_{\Theta(B+1)}(T_i)
 }
 { g(T_i)    \sum_{j=1}^{N} \frac{\exp\left[-\lambda r_{t}(T_j)
\right]\mathbf{1}_{\Theta(B+1)}(T_j)}{ g(T_j)} }
\xrightarrow[N\rightarrow\infty]{a.s.}
\hat{\theta}^{\tau}_{B,\lambda}[t].
 $$
Remark that this is particularly convenient as we only simulate the sample
$T_1$, ..., $T_N$ once and
we can use the previous formula to approximate $\hat{\theta}^{\tau}_{B,\lambda}[t]$
for several values
of $\tau$.

\subsection{Results}

\begin{figure}[!!h]
\begin{center}
\centering
\includegraphics*[width=5cm]{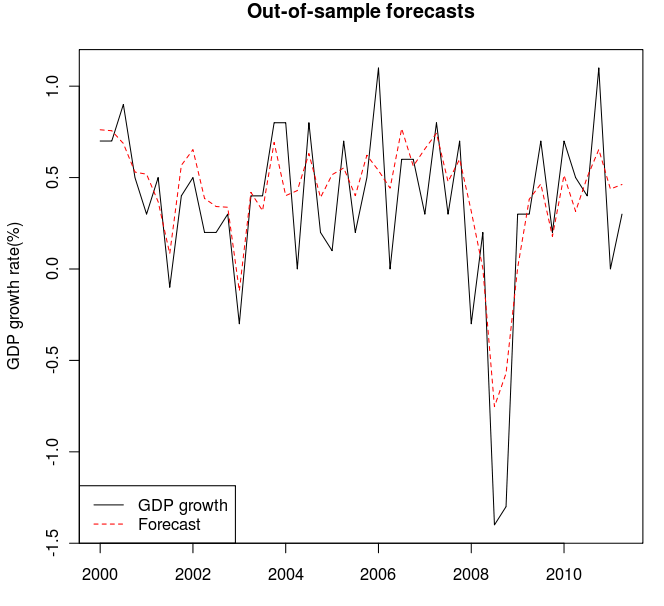}
\caption{\label{fig05} \textsf{French GDP online prediction using the quantile
loss function with $\tau=0.5$.}}
\end{center}
\end{figure}

\begin{figure}[!!h]
\begin{center}
\centering
\begin{tabular}{c c} 
\includegraphics*[width=5cm]{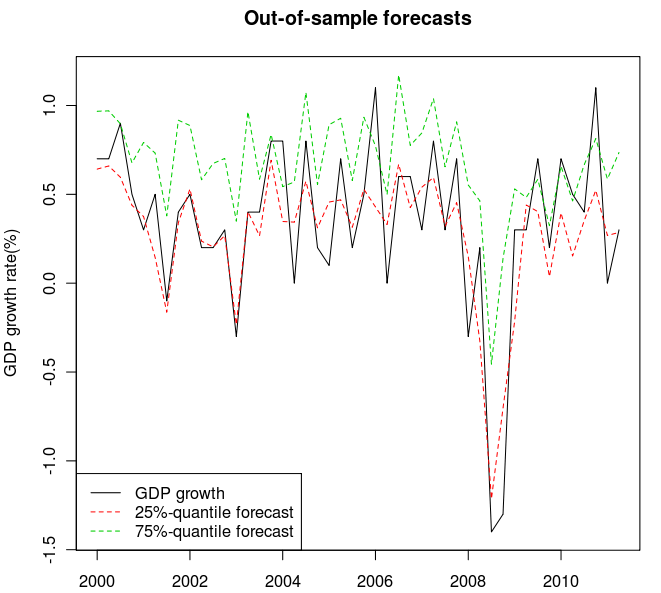} & \includegraphics*[width=5cm]{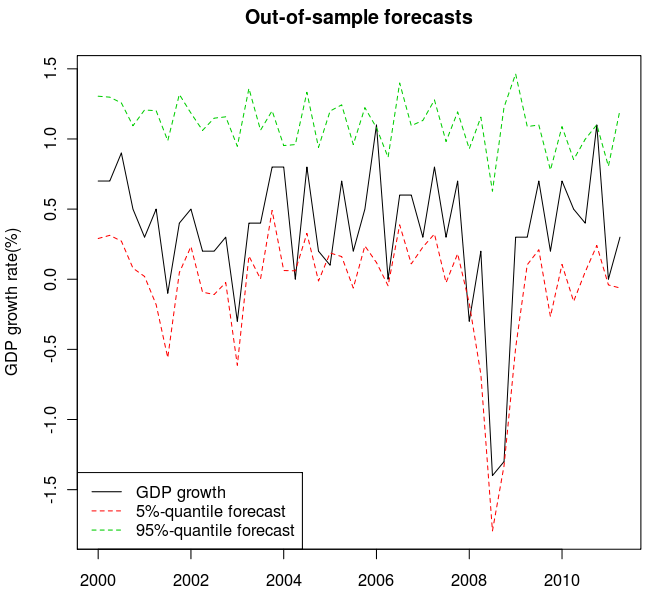}
\end{tabular}
\caption{\label{fig025} \textsf{French GDP online $50\%$-confidence intervals (left)
and $90\%$-confidence intervals (right).}}
\end{center}
\end{figure}

The results are shown in
Figure \ref{fig05} for prediction, $\tau=0.5$, in Figure \ref{fig025} for
confidence interval of order $50\%$, i.e. $\tau=0.25$ and $\tau=0.75$ (left) and for
confidence interval of order $90\%$, i.e. $\tau=0.05$ and $\tau=0.95$ (right).
We report only the results for the period 2000-Q1 to 2011-Q3 (using the period
1988-Q1 to 1999-Q4 for learning).

Note that we can compare the ability of our predictor $\widehat{\theta}_{B,\lambda}^{0.5}$
with the predictor used in \cite{stage} that relies on a least square estimation
of~\eqref{modele}, that we will denote by $\hat{\theta}^*$. Interestingly, both are
quite similar but $\widehat{\theta}_{B,\lambda}^{0.5}$ is a bit more precise. We remind
that
$$
\begin{array}{c c}
\text{mean abs. pred. error} &=  \frac{1}{n}\sum_{t=1}^{n} \left| \Delta GDP_{t} -
f_{\hat{\theta}^{0.5}_{B,\lambda(t)}[t]}(X_{t-1},X_{t-2}) \right|
\\
\text{mean quad. pred. error} &= \frac{1}{n}\sum_{t=1}^{n} \left[ \Delta GDP_{t} -
f_{\hat{\theta}^{0.5}_{B,\lambda(t)}[t]}(X_{t-1},X_{t-2}) \right]^2.
\end{array}
$$

\begin{center}
\begin{tabular}{|c|c|c|}
\hline
 Predictor & Mean absolute prevision error & Mean quadratic prevision error \\
$\widehat{\theta}_{B,\lambda}^{0.5}$ & $0.22360$ & $0.08033$ \\
$\widehat{\theta}_{\star}$ & $0.24174$  & $0.08178 $ \\
\hline
\end{tabular}
\end{center}

We also report the frequency of realizations of the GDP falling above the
predicted $\tau$-quantile for each $\tau$.
Note that this quantity should be close to $\tau$.

\begin{center}
\begin{tabular}{|c|c|}
\hline 
 Estimator & Frequency \\
$\widehat{\theta}_{B,\lambda}^{0.05}$ & $0.065$ \\
$\widehat{\theta}_{B,\lambda}^{0.25}$ & $0.434$ \\
$\widehat{\theta}_{B,\lambda}^{0.5}$  & $0.608$ \\
$\widehat{\theta}_{B,\lambda}^{0.75}$ & $0.848$ \\
$\widehat{\theta}_{B,\lambda}^{0.95}$ & $0.978$ \\
\hline
\end{tabular}
\end{center}

As the INSEE did, we
miss the value of the 2008 crisis. However, it is interesting to note that our confidence
interval shows that our prediction at this date is less reliable than the previous ones:
so, at this time, the forecaster could have been aware of some problems in their predictions.

\section{Conclusion}

We proposed some theoretical results to extend learning theory to the context of weakly
dependent time series. The method showed good results on an application to GDP forecasting.
A next step will be to give theoretical results on the online risk of our method,
e.g. using tools from \cite{Catoni2004,gerchi}. From both theoretical and practical
perspective, an adaptation with respect to the dependence coefficient $\theta_{\infty,n}(1)$
would also be really interesting but is probably a more difficult objective.

\section*{Acknowledgments}

We deeply thank Matthieu Cornec (INSEE) who
provided the data and the model in the application. We also thank
Pr. Paul Doukhan (Univ. Cergy-Pontoise) for his helpful comments.
Research partially supported by the ``Agence
Nationale pour la Recherche'', grant ANR-09-BLAN-0128 ``PARCIMONIE''.

\section{Proofs}
\label{sectionproofs}

\subsection{Some preliminary lemmas}

Our main tool is Rio's Hoeffding type inequality \cite{Rio2000a}.
The reference \cite{Rio2000a}
is written in French and unfortunately, up to our knowledge, there is no
English version of this result. So for the sake of completeness, we
provide this result.
\begin{lemma}[Rio \cite{Rio2000a}]
\label{RIO}
Let $h$ be a function $(\mathbb{R}^{p})^{n}\rightarrow\mathbb{R}$ such that for all $x_1$,
..., $x_n$, $y_1$, ..., $y_n\in\mathbb{R}^{p}$,
\begin{equation}
\label{lipcondrio}
         |h(x_1,\ldots,x_n)- h(y_1,\ldots,y_n)|
\leq \sum_{i=1}^{n} \|x_i-y_i\|.
\end{equation}
Then for any $t>0$ we have
$$
\mathbb{E}\left(e^{t\left\{\mathbb{E}\left[h(X_{1},\ldots,X_{n})\right]
 - h(X_{1},\ldots,X_{n}) \right\}}\right)
\leq e^{\frac{t^{2} n \left(\mathcal{B} + \theta_{\infty,n}(1)\right)^{2}}{2}}
.
$$
\end{lemma}
Others exponential inequalities can be used
to obtain PAC-Bounds in the context of time series: the inequalities in \cite{Doukhan1994}
for mixing time series, and \cite{Dedecker2007a,devdep} under weakest ``weak dependence''
assumptions, \cite{Seldin} for martingales. However, Lemma~\ref{RIO} is
particularly convenient here, and will lead to optimal learning rates.

\begin{lemma}
\label{XIAOYIN}
We remind that $\kappa = K(1+L) (\mathcal{B} + \mathcal{C})/\sqrt{2}$.
We assume that Assumptions {\bf LipLoss}, {\bf WeakDep} and {\bf Lip} are
satisfied.
For any $\lambda>0$ and $\theta\in\Theta$,
$$ \mathbb{E}\left( e^{\lambda [R(\theta) - r_n (\theta)] } \right) \leq
e^{\frac{\lambda^{2} \kappa^{2}}{n\left(1-\frac{k}{n}\right)^2}}
 \text{  and  }
\mathbb{E}\left( e^{\lambda [r_n(\theta) - R (\theta)] } \right) \leq
e^{\frac{\lambda^{2} \kappa^{2}}{n\left(1-\frac{k}{n}\right)^2}} .$$
\end{lemma}
\noindent {\it Proof of Lemma~\ref{XIAOYIN}.}
Let us fix $\lambda>0$ and $\theta\in\Theta$. Let us define the function $h$ by:
$$
h(x_1, \ldots , x_n)=\frac{1}{K(1+L)}\sum_{i=k+1}^n
\ell(f_{\theta}(x_{i-1},\ldots,
x_{i-k}),x_i).
$$
We now check that $h$ satisfies~\eqref{lipcondrio}, remember that $\ell(x,x')=
g(x-x')$ so
\begin{align*}
 & \Bigl|  h\left(  x_{1},\ldots, x_{n}\right)  - h\left(y_{1},\ldots
y_{n}\right)\Bigr| \\
 & \quad \leq \frac{1}{K(1+L)}\sum_{i=k+1}^{n} \Bigl|
g(f_{\theta}(x_{i-1},\ldots,x_{i-k})-x_{i})
                                 -g(f_{\theta}(y_{i-1},\ldots,y_{i-k})-y_{i})
\Bigr| \\ 
 & \quad \leq \frac{1}{1+L}\sum_{i=k+1}^{n} \Bigl\|
\bigl(f_{\theta}(x_{i-1},\ldots,x_{i-k})-x_{i}\bigr)                                
-\bigl(f_{\theta}(y_{i-1},\ldots,y_{i-k})-y_{i}\bigr) \Bigr\|
\end{align*}
where we used Assumption {\bf LipLoss} for the last inequality. So we have
\begin{align*}
 & \Bigl|  h\left(  x_{1},\ldots, x_{n}\right)  - h\left(y_{1},\ldots
y_{n}\right)\Bigr| \\
 & \quad \leq \frac{1}{1+L}\sum_{i=k+1}^{n} \biggl( \Bigr\|
f_{\theta}(x_{i-1},\ldots,x_{i-k})
                        - f_{\theta}(y_{i-1},\ldots,y_{i-k}) \Bigr\|  + \Bigl\|
x_i - y_i \Bigr\| \biggr) \\
 & \quad \leq \frac{1}{1+L}\sum_{i=k+1}^{n} \left( \sum_{j=1}^{k} a_{j}(\theta)
                                        \| x_{i-j} - y_{i-j}\|
                 +  \| x_i - y_i \| \right) \\
 & \quad \leq \frac{1}{1+L}\sum_{i=1}^{n} \left(1 +
             \sum_{j=1}^{k}a_{j}(\theta) \right) \| x_i - y_i \| 
 \leq \sum_{i=1}^{n} \|x_i-y_i \|
 \end{align*}
where we used Assumption {\bf Lip}. So we can apply Lemma~\ref{RIO} with
$h(X_1,\ldots,X_n)=\frac{n-k}{K(1+L)}
r_{n}(\theta)$, $\mathbb{E}(h(X_1,\ldots,X_n))=\frac{n-k}{K(1+L)} R(\theta)$, and
$t=K(1+L)\lambda/(n-k)$:
\begin{equation*}
\mathbb{E}\left(e^{\lambda\left[R(\theta)
 - r_{n}(\theta) \right]}\right)
\leq e^{\frac{\lambda^{2} K^{2} (1+L)^{2}
    \left(\mathcal{B} + \theta_{\infty,n}(1)\right)^{2}}{2 n
    \left(1-\frac{k}{n}\right)^{2} }}
\leq e^{\frac{\lambda^{2}K^{2} (1+L)^{2}
    \left(\mathcal{B} + \mathcal{C}\right)^{2}}{2 n
    \left(1-\frac{k}{n}\right)^{2} }}
\end{equation*}
by Assumption {\bf WeakDep}. This ends the proof of the first
inequality.
The reverse inequality is obtained by replacing the function $h$ by $-h$.
\hfill $\blacksquare$

We also remind the following result concerning the Kullback divergence.
\begin{lemma}
\label{LEGENDRE}
For any $\pi\in\mathcal{M}_{+}^{1}(E)$,
for any measurable upper-bounded function $h:E\rightarrow\mathbb{R}$ we have:
\begin{equation} \label{lemmacatoni}
\pi[\exp
(h)]=\exp\left(\sup_{\rho\in\mathcal{M}_{+}^{1}(E)}\biggl(\rho
[h]-\mathcal{K}(\rho,\pi)\biggr)\right).
\end{equation}
Moreover, the supremum with respect to $\rho$ in the right-hand side is
reached for the Gibbs measure
$\pi\{h\}$ defined by $ \pi\{h\}({\rm d}x) =  e^{h(x)} \pi({\rm d}x) / \pi[\exp(h)]$.
\end{lemma}
Actually, it seems that in the case of discrete probabilities, this
result was already known by Kullback (Problem 8.28 of Chapter 2 in
\cite{kullback}).
For a complete proof in the general case, we refer the reader for example
to \cite{Catoni2003,Catoni2007}.
We are now ready to state the following key result.
\begin{lemma}
\label{PACBAYES}
 Let us assume that Assumptions {\bf LipLoss}, {\bf WeakDep} and {\bf Lip} are
satisfied.
 Let us fix $\lambda>0$. Let $k$ be defined as in Lemma~\ref{XIAOYIN}.
 Then,
\begin{equation}
\label{eqPACBAYES}
\mathbb{P} \left\{
\begin{array}{l}
\forall \rho \in\mathcal{M}_{+}^{1}(\Theta), \\
 \int R {\rm d} \rho \leq \int r_n {\rm d}\rho
       + \frac{\lambda \kappa^2 }{n \left(1-\frac{k}{n}\right)^{2}}
     + \frac{\mathcal{K}(\rho,\pi) +
\log\left(\frac{2}{\varepsilon}\right)}{\lambda}
\\
\text{ and }
\\
 \int r_n {\rm d} \rho \leq \int R {\rm d}\rho
       + \frac{\lambda \kappa^2 }{n \left(1-\frac{k}{n}\right)^{2}}
     + \frac{\mathcal{K}(\rho,\pi) +
\log\left(\frac{2}{\varepsilon}\right)}{\lambda}
\end{array}
\right\}
\geq 1-\varepsilon.
\end{equation}
\end{lemma}
\noindent {\it Proof of Lemma~\ref{PACBAYES}.}
 Let us fix $\theta>0$ and $\lambda>0$, and apply the first inequality of
Lemma~\ref{XIAOYIN}. We have:
$$
\mathbb{E}e^{\lambda \left[R(\theta) - r_n (\theta) - \frac{\lambda
\kappa^{2}}{n\left(1-\frac{k}{n}\right)^2}\right] }
          \leq 1,
$$
and we multiply this result by $\varepsilon/2$ and integrate it with
respect to $\pi({\rm d}\theta)$. Fubini's Theorem gives:
$$
\mathbb{E}\int e^{\lambda \left[R(\theta) - r_n (\theta)\right]
- \frac{\lambda^2 \kappa^{2}}{n\left(1-\frac{k}{n}\right)^2} - \log
\left(\frac{2}{\varepsilon}\right)}
       \pi({\rm d}\theta) \leq \frac{\varepsilon}{2}.
$$
We apply Lemma~\ref{LEGENDRE} and we get:
$$
\mathbb{E} e^{\sup_{\rho} \left\{ \lambda \int \left[R(\theta) - r_n
(\theta)\right] \rho({\rm d}\theta)
 - \frac{\lambda^2 \kappa^{2}}{n\left(1-\frac{k}{n}\right)^2} - \log
\left(\frac{2}{\varepsilon}\right) - \mathcal{K}(\rho,\pi)
 \right\}  }  \leq \frac{\varepsilon}{2}.
$$
As $e^{x}\geq \mathbf{1}_{\mathbb{R}_{+}}(x)$, we have:
$$
\mathbb{P} \left\{
\sup_{\rho} \left\{ \lambda \int \left[R(\theta) - r_n (\theta)\right] \rho({\rm
d}\theta)
 - \frac{\lambda^2 \kappa^{2}}{n\left(1-\frac{k}{n}\right)^2} - \log
\left(\frac{2}{\varepsilon}\right) - \mathcal{K}(\rho,\pi)
 \right\} \geq 0
 \right\} \leq \frac{\varepsilon}{2}.
$$
Now, we follow the same proof again but starting with the second inequality of
Lemma~\ref{XIAOYIN}. We obtain:
$$
\mathbb{P} \left\{
\sup_{\rho} \left\{ \lambda \int \left[r_n(\theta) - R (\theta)\right] \rho({\rm
d}\theta)
 - \frac{\lambda^2 \kappa^{2}}{n\left(1-\frac{k}{n}\right)^2} - \log
\left(\frac{2}{\varepsilon}\right) - \mathcal{K}(\rho,\pi)
 \right\} \geq 0
 \right\} \leq \frac{\varepsilon}{2}.
$$
A union bound ends the proof.
\hfill $\blacksquare$

\subsection{Proof of Theorems \ref{main_result} and \ref{thm_appli}}

\noindent {\it Proof of Theorem~\ref{main_result}.}
We apply Lemma~\ref{PACBAYES}. So, with probability at least $1-\varepsilon$
we are on the event given by~\eqref{eqPACBAYES}. From now, we work on that
event. The first inequality of~\eqref{eqPACBAYES}, when applied to $
\hat{\rho}_{\lambda}({\rm d}\theta)$, gives
$$
 \int R(\theta) \hat{\rho}_{\lambda}({\rm d}\theta) \leq  \int r_n (\theta)
\hat{\rho}_{\lambda}({\rm d}\theta)
 + \frac{\lambda \kappa^{2}}{n\left(1-\frac{k}{n}\right)^2} +\frac{1}{\lambda}
\log \left(\frac{2}{\varepsilon}\right)
+ \frac{1}{\lambda} \mathcal{K}(\hat{\rho}_{\lambda},\pi).
$$
According to Lemma \ref{LEGENDRE} we have:
$$
\int r_n (\theta) \hat{\rho}_{\lambda}({\rm d}\theta)
 + \frac{1}{\lambda} \mathcal{K}(\hat{\rho}_{\lambda},\pi) = \inf_{\rho}
\left(\int r_n (\theta) \rho({\rm d}\theta)
+ \frac{1}{\lambda}\mathcal{K}(\rho,\pi) \right)
$$
so we obtain
\begin{equation} \label{eq1}
 \int R(\theta) \hat{\rho}_{\lambda}({\rm d}\theta)  \leq \inf_{\rho} \left[
\int r_n (\theta) \rho({\rm d}\theta)
 + \frac{\lambda \kappa^{2}}{n\left(1-\frac{k}{n}\right)^2} +
\frac{\mathcal{K}(\rho,\pi) + \log\left(\frac{2}{\varepsilon}\right)}{\lambda}
 \right].
\end{equation}
We now want to bound from above $r(\theta)$ by $R(\theta)$. 
Applying the second inequality of~\eqref{eqPACBAYES} and plugging it into
Inequality~\ref{eq1}
gives
$$
    \int R(\theta) \hat{\rho}_{\lambda}({\rm d}\theta)  \leq \inf_{\rho} \left[
\int R {\rm d}\rho+\frac{2}{\lambda}\mathcal{K}(\rho,\pi)
 + \frac{2\lambda \kappa^{2}}{n\left(1-\frac{k}{n}\right)^2} + \frac{2}{\lambda}
\log\left(\frac{2}{\varepsilon}\right) \right].
$$
We end the proof by the remark that $\theta\mapsto R(\theta)$ is convex and
so by Jensen's inequality
$ \int R(\theta) \hat{\rho}_{\lambda}({\rm d}\theta) \geq
R\left( \int \theta \hat{\rho}_{\lambda}({\rm d}\theta) \right)
= R(\hat{\theta}_{\lambda}) .$
\hfill $\blacksquare$

\noindent {\it Proof of Theorem~\ref{thm_appli}.}
We can apply Theorem \ref{main_result} with $R=R_\tau$. 
We have, with probability at least $1-\varepsilon$,
\begin{equation*}
  R^{\tau} (\hat{\theta}_{B,\lambda}^{\tau}) \leq
\inf_{\rho\in\mathcal{M}_{+}^{1}(\Theta)} \left[ \int R^{\tau} {\rm d}\rho
 + \frac{2\lambda \kappa^{2}}{n\left(1-\frac{2}{n}\right)^2} +
\frac{2\mathcal{K}(\rho,\pi) +  2 \log\left(\frac{2}{\varepsilon}\right)
}{\lambda} \right].
\end{equation*}
Now, let us fix $\delta\in(0,1]$ and $\theta\in\Theta(B)$.
We define the probability distribution $\rho_{\theta,\delta}$ as the uniform
probability measure on the set $ \{T \in \mathbb{R}^{4} ,\quad \|\theta-T\|_{1}
\leq \delta \}$.
Note that $\rho_{\theta,\delta} \ll \pi_{B} $ as $\pi_{B}$ is defined as uniform
on $\Theta(B+1)\supset \Theta(B+\delta)$.
Then:
\begin{equation}
\label{stepproof1}
 R^{\tau} (\hat{\theta}_{B,\lambda}^{\tau}) \leq \inf_{
\tiny{
\begin{array}{c}
\theta\in\Theta(B) \\
\delta>0
\end{array}
}
} \left[ \int R^\tau {\rm d}\rho_{\theta,\delta}
 + \frac{2\lambda \kappa^{2}}{n\left(1-\frac{2}{n}\right)^2} +
\frac{2\mathcal{K}(\rho_{\theta,\delta},\pi)
       +  2 \log\left(\frac{2}{\varepsilon}\right)
}{\lambda} \right].
\end{equation}
Now, we have to compute or to upper-bound all the terms in the right-hand side
of this inequality.
First, note that:
\begin{multline}
\label{stepproof2}
\int R^\tau {\rm d}\rho_{\theta,\delta}
 = \int_{ \{\|\theta-T\|_{1} \leq \delta \} }
       R^\tau (T) {\rm d}\rho_{\theta,\delta}(T) 
\\
\leq R^{\tau}(\theta) + 2 \mathcal{B} \delta \max(\tau,1-\tau)\leq
R^{\tau}(\theta) + 2 \mathcal{B} \delta .
\end{multline}
We compute $
 \mathcal{K}(\rho_{\theta,\delta},\pi_{B}) = 3
\log\left(\frac{B+1}{\delta}\right)
$ and plug this with \eqref{stepproof2}
 into \eqref{stepproof1} to
obtain:
\begin{equation*}
 R^{\tau} (\hat{\theta}_{B,\lambda}^{\tau}) 
\leq \inf_{\theta,\delta}\left\{ R^{\tau}(\theta)
 + 2 \left[\frac{\lambda \kappa^{2}}{n\left(1-\frac{2}{n}\right)^2} +
\mathcal{B} \delta + \frac{3\log \left(\frac{B+1}{\delta}\right)+
\log\left(\frac{2}{\varepsilon}
\right)}{\lambda} \right] \right\}.
\end{equation*}
It can easily be seen that the minimum of the right-hand side w.r.t. $\delta$ is
reached for $\delta=3/(\mathcal{B}\lambda)\leq 1$ as soon as $\lambda$ is large enough,
and so:
\begin{equation*}
 R^{\tau} (\hat{\theta}_{B,\lambda}^{\tau}) 
\leq \inf_{\theta}\left\{ R^{\tau}(\theta) + \frac{2\lambda
\kappa^{2}}{n\left(1-\frac{2}{n}\right)^2}
  + \frac{6\log \left(\frac{(B+1)\mathcal{B}\lambda {\rm e}}{3}\right)+ 2
\log\left(\frac{2}{\varepsilon}
\right)}{\lambda}\right\}.
\end{equation*}
We finally propose
$\lambda = \sqrt{3n}/\kappa$, this leads to:
\begin{equation*}
 R^{\tau} (\hat{\theta}_{B,\lambda}^{\tau}) 
\leq \inf_{\theta}\left\{ R^{\tau}(\theta) +
\frac{2\sqrt{3} \kappa}{\sqrt{n}} \left[
\frac{1}{\left(1-\frac{2}{n}\right)^2}
  + \log \left(\frac{(B+1)\mathcal{B} {\rm e}}{\kappa} \sqrt{\frac{n}{3}}
\right)
+ \frac{\log\left(\frac{2}{\varepsilon}\right)}{3}
   \right]\right\}.
\end{equation*}
Remark that the condition $\delta<1$ is satisfied as soon as $n>\kappa^2 /
(3\mathcal{B}^{2})$.
Moreover, $ \forall n\geq 10, \quad 1/\left(1-\frac{2}{n}\right)^2
\leq \frac{25}{16} $
and we can re-organize the terms to obtain:
\begin{equation*}
 R^{\tau} (\hat{\theta}_{B,\lambda}^{\tau}) 
\leq \inf_{\theta}\left\{ R^{\tau}(\theta) +
\frac{2\sqrt{3} \kappa}{\sqrt{n}} \left[
2.25 + \log\left(\frac{(B+1)\mathcal{B} \sqrt{n}}{\kappa}\right)
+ \frac{\log\left(\frac{1}{\varepsilon}\right)}{3}
   \right]\right\}.
\end{equation*}
\hfill $\blacksquare$

\bibliographystyle{splncs}
\bibliography{biblio}

\end{document}